\documentclass[12pt]{article}
\usepackage{amsmath,amsthm,amsfonts,amssymb}
\def\Diff{{\hbox{\rm Diff}}}

\def\Hom{{\hbox{Hom}}}
\def\id{{\rm id}}

\def\a{\alpha}

\def\epsilon{\varepsilon}

\def\la{\lambda}

\def\phi{\varphi}

\newtheorem{theorem}{Theorem}[section]

\def\A{{\mathcal A}}
\def\B{{\mathcal B}}
\def\Z{{\mathbb Z}}

\def\R{{\mathbb R}}

\title{{\bf Classification of operator algebraic conformal field theories}}

\author{{\sc Yasuyuki Kawahigashi}
\footnote{The author was supported in part by JSPS Grant.}\\
Department of Mathematical Sciences\\
University of Tokyo, Komaba, Tokyo, 153-8914, Japan\\
e-mail: {\tt yasuyuki@ms.u-tokyo.ac.jp}}
\begin{document}
\date{}
\maketitle

\begin{abstract}
We give an exposition on the current status of classification of
operator algebraic conformal field theories.  We explain roles
of complete rationality and $\alpha$-induction for nets of subfactors
in such a classification and present the current classification
result, a joint work with R. Longo,
for the case with central charge less than $1$, where
we have a complete classification list consisting of the
Virasoro nets, their simple current extensions of index 2, and four 
exceptionals.  Two of the four exceptionals appear to be new.
\end{abstract}

\section{Introduction}

Classification of conformal field theory is obviously one of the
most important and exciting problems in mathematical physics today.
Here we review the current status of the classification theory of
1-dimensional operator algebraic conformal field theories.

First, we explain what we mean by 
``1-dimensional operator algebraic conformal field theories'' and 
how they naturally appear in the setting of $1+1$-dimensional
conformal field theory.
We are now in the framework of algebraic quantum field theory in the
sense of \cite{H}.  In the algebraic quantum field theory, we 
consider a net of von Neumann algebras $\{\B(O)\}$ where $O$ is
a certain bounded region in the Minkowski space, and the physical
idea is that the von Neumann
algebra $\B(O)$ is generated by physical quantities observable
in the region $O$.  From a physical viewpoint,
the most natural setting is that we have a 4-dimensional Minkowski
space and $O$ is a double cone, which is defined to be a set of the
form $(x + V_+)\cap(y + V_-)$, where
$$V_\pm = \{ z=(z_0,z_1,z_2,z_3)\in \R^4 \mid
z_0^2 - z_1^2 - z_2^2 - z_3^2 > 0,  \pm z_0 > 0\}.$$
We say a {\sl net} because the index set of double cones
is directed with respect to inclusions.
Then we impose some physically natural set of axioms
on the net $\B$, as explained below.
So mathematically speaking, this is a study
of a family of von Neumann algebras (on a fixed Hilbert space)
parametrized by double cones subject to certain set of axioms.
Although the 4-dimensional Minkowski space seems most natural
at first sight, it has been known that a net on
the 2-dimensional Minkowski space produces interesting structures
both from physical and mathematical viewpoints.
Such a net is called a two-dimensional net.
Using a new set of coordinates $(x_0+x_1, x_0-x_1)$
in the two-dimensional Minkowski space, a double cone $O$ becomes
a product of two intervals $I=(a,b)$, $J=(c,d)$.  We now start with
a net $\B$ of von Neumann algebras on the two-dimensional
Minkowski space.  Then as in \cite{R2}, we obtain a decomposition
$\A_L(I)\otimes \A_R(J)\subset \B(I\times J)$, where
$\A_L$, $\A_R$ are nets of von Neumann algebras on one-dimensional
spaces.  Such one-dimensional nets are called {\sl chiral theories}.
In this way, we are naturally led to a study of one-dimensional
nets of von Neumann algebras, starting from a study of two-dimensional
nets of von Neumann algebras.  As a part of the set of axioms for
algebraic quantum field theory, we assume that we have a certain
group representing  space-time symmetry and its projective
unitary representation on the Hilbert space on which our net
of von Neumann algebras act.  We require compatibility of this
projective unitary representation and the net of von Neumann algebras
in an appropriate sense and this compatibility condition is
called {\sl covariance}.  A space-time symmetry of the
two-dimensional Minkowski space passes to a symmetry of
the one-dimensional space in the above decomposition, so we now
consider a symmetry group of the one-dimensional space.  It is
often more convenient to compactify the one-dimensional space
and use $S^1$ rather than $\R$, and we do so here.
There can be more than one choice for such a 
symmetry group of $S^1$, and we often use the name ``conformal
field theory'' when we choose the M\"obius group $PSL(2,\R)$ as the
symmetry group.
In this paper, we use this M\"obius group, or even more
strongly, the diffeomorphism group $\Diff(S^1)$ later.  (See \cite{KL}
for a precise definition of the diffeomorphism covariance.)
Thus our set of axioms for one-dimensional net of von Neumann
algebras is now described as follows.  

By an {\sl interval}, we mean   a non-empty, non-dense,
open, and connected set in $S^1$.  We study a family of von Neumann
algebras $\A(I)$ on a fixed Hilbert space $H$ parameterized by
intervals $I$.  (Now the set of intervals is not directed with
respect to inclusions, so the
terminology ``net'' is not appropriate.  The name ``precosheaf''
is also often used instead, because of this reason, but the name
``net'' has been used in many cases, and we also use it here.)
Our set of axioms is as follows.  (See \cite{GL2,KL}
for precise descriptions of these properties.)
\begin{enumerate}
\item Isotony
\item Locality
\item Covariance
\item Positivity of the energy
\item Unique existence of a vacuum vector
\end{enumerate}
Such a family of von Neumann algebras is also called a
{\sl local conformal net}.  Then the {\sl Haag duality},
$\A(I')=\A(I)'$, automatically holds, where the prime on the left hand
side means the interior of the complement of the interval
$I$, and the prime on the right hand side denotes the
commutant of $\A(I)$.

The main object of mathematical studies of such nets is their
classification.  By a ``classification'', we mean obtaining
a complete invariant up to isomorphism and a complete list of
all the isomorphism classes.  Both problems are very difficult, but
we explain the first step to these problems here.  The invariant
to the classification problem we now use is a representation theory of 
a net.  Under the above set of axioms and the split
property, which is often assumed as an additional axiom,
each von Neumann algebra $\A(I)$ is an injective type III$_1$ factor.  
In general, all representations of a type III factor on  separable
Hilbert spaces are unitarily equivalent, so representation
theory for a single algebra $\A(I)$ is trivial.  However,
we obtain non-trivial and useful information by considering
representations of a {\sl family} of type III factors.
That is, we consider a family of representations $\pi_I$ of
$\A(I)$ on one Hilbert space parameterized by intervals $I$ on the circle,
with a natural compatibility condition and
a covariance property with respect to the symmetry group.
(By \cite{GL1}, a covariance property with respect to
the M\"obius group automatically holds, if we have strong
additivity and the statistical dimension, which is explained below,
of the representation is finite.)  By adapting the
standard Doplicher-Haag-Roberts (DHR) theory \cite{DHR} to our
one-dimensional setting, we get a category of representations and
each representation is realized as a DHR endomorphism.
Then we have a notion of irreducibility, direct sums,
irreducible decomposition, a statistical dimension, conjugates,
unitary equivalence and a tensor product.  It behaves like a
category of unitary representations of a compact group.
Here a statistical dimension plays a role
of a usual dimension, but the values are now in the set of
positive real numbers.  A tensor product is given by a
composition of DHR endomorphisms.  In this way, we obtain a
strict $C^*$-tensor category with conjugates, subobjects, and
direct sums in the sense of \cite{DR,LRo}.

In operator algebra theory, we often have an algebraic
invariant from a certain representation theoretic consideration,
and it becomes a complete invariant under an extra assumption
of {\sl amenability} in some sense, though proving completeness
of such an invariant requires very deep analytic arguments.
Also in this setting of algebraic quantum field theory, we expect
that a representation category, possibly with some additional
related data, should give a complete invariant of local conformal 
nets with some kind of amenability assumption, but results of such 
a general type have not been obtained and seem very difficult to prove.  
In this paper, we explain the first classification result \cite{KL}
in this direction under an extra assumption on a central charge,
a numerical invariant of a net that is also ``representation
theoretic''.

When we assume covariance with respect to $\Diff(S^1)$ rather
than $PSL(2,\R)$ for a local conformal net $\A$,
we can define a central charge $c$ for $\A$ as follows.
The corresponding infinite dimensional Lie algebra of the 
infinite dimensional Lie group $\Diff(S^1)$ is the
celebrated Virasoro algebra, which is generated by
$\{L_m\}_{m\in\Z}$ and a central element $c$ subject
to the following relations.
$$[L_m, L_n]=(m-n)L_{m+n}+\frac{c}{12}(m^3-m)\delta_{m,-n},$$
where $m,n\in\Z$.  In our setting, we obtain a
representation of the Virasoro algebra from a local
conformal net with diffeomorphism covariance,
and we can define the central charge $c$, which is a numerical
value.  We then
define this value to be the central charge of the net $\A$.
It is clearly a numerical invariant of $\A$.  It has been shown by
Friedan-Qiu-Shenker
\cite{FQS} that this central charge value is in
$$\{1-6/m(m+1)\mid m=2,3,4,\dots\}\cup[1,\infty)$$ and
the values $\{1-6/m(m+1)\mid m=2,3,4,\dots\}$ have been 
realized by Goddard-Kent-Olive
\cite{GKO}.  (The values in $[1,\infty)$ are easier to
realize.)  Similarity of this restriction of
possible values to the restriction of the values of the
Jones indices of subfactors in \cite{J} has been pointed by
Jones himself since early days of subfactor theory.  Furthermore,
in subfactor theory,
we have classification of {\sl paragroups} of index less than 4
with the Dynkin diagrams
$A_n$, $D_{2n}$, $E_6$, $E_8$ by Ocneanu \cite{O1}, which
produces a classification of subfactors of the hyperfinite II$_1$
factor with index less than 4 by Popa's analytic classification
result \cite{P}.  (See the book \cite{EK} for details on this
topic.)  Also in the theory of Virasoro algebra, we
have a classification involving $A$-$D$-$E$ Dynkin diagrams,
that is, classification of {\sl modular invariants} by
Cappelli-Itzykson-Zuber \cite{CIZ}.  (See also the book
\cite{DMS} for more on modular invariants.)  The similarity between
the two classification is lised in Table \ref{tab-Vir}.

\begin{table}[htbp]
\begin{center}
\begin{tabular}{|c|c|}\hline
Virasoro algebra & Subfactors \\ \hline
central charge $c$ & Jones index \\ 
$c<1\Rightarrow c=1-6/m(m+1)$ 
\cite{FQS} &  ${\rm index}<4 \Rightarrow {\rm index}=4\cos^2\pi/m$
\cite{J} \\ \hline
classification of modular invariants
\cite{CIZ} & classification of paragroups \cite{O1}  \\ 
with pairs of $A$-$D$-$E$ diagrams & with $A$-$D$-$E$ diagrams \\ \hline
\end{tabular}
\caption{Similarity between Virasoro algebra and subfactors}
\label{tab-Vir}
\end{center}
\end{table}

However, a real relation beyond formal similarity of the two
classification theories was unknown.  Furthermore, even on a 
formal level of similarity, it is not clear at all why we have 
pairs of Dynkin diagrams for Virasoro algebra and single
Dynkin diagrams for subfactors.  Our results in \cite{KL}
on classification of local conformal nets, with diffeomorphism
covariance and central charge
less than 1, give unification of these two classification
theories.

We also remark a relation between such a study of local conformal nets
and $E_0$-semigroups.  The $E_0$-semigroups appearing here are not
those of $B(H)$, but of the injective type III$_1$-factor.
Consider an inclusion $N\subset M$ of von Neumman algebra with
a common cyclic separating vector $\Omega$.  We have a modular
automorphism group $\sigma^M_t$ with respect to the pair
$(M,\Omega)$.  Wiesbrock studied ``half-sided modular inclusions''
where we have $\sigma^M_t(N)\subset N$ for $t\leq 0$.  This gives an
$E_0$-semigroup of $N$, which is typically an injective type 
III$_1$-factor in the setting appearing from a local conformal
net.  Wiesbrock \cite{Wi} and Guido-Longo-Wiesbrock 
\cite[Theorem 1.2]{GLW} have shown that
we have a triple of half-sided modular inclusions from a local
conformal net on the circle, by splitting the circle into three
intervals after removing three distinct points, and we can also
recover the local conformal net from such a triple.  In this sense,
study of local conformal nets on the circle is ``translated'' to that
of special type of $E_0$-semigroups of the injective type
III$_1$-factor.

\section{Completely rational nets and $\a$-induction}

We have mentioned that the invariant of a local conformal
net we study is its representation theory, but in general, 
it is difficult to control representation theory of a 
net.  A set of additional axioms has been
proposed by us \cite{KLM} in order to obtain a ``good''
representation theory as follows.  (See \cite{KLM} for the
precise definitions.)
\begin{enumerate}
\item Split property
\item Strong additivity
\item Finite index condition
\end{enumerate}
The first two conditions are rather well-known in algebraic
quantum field theory.
The definition of the new third condition, finite index condition,
is as follows.  Remove four distinct points from the circle $S^1$
and label the remaining four connected components as
$I_1, I_2, I_3, I_4$ in a counterclockwise order.
By locality, we have an inclusion
$\A(I_1)\vee \A(I_3)\subset (\A(I_2)\vee\A(I_4))'$ of von 
Neumann algebras, which are now automatically injective
type III-factors.  The finite index condition means 
that the index of this subfactor is finite.  (This condition
is independent of choices of the four intervals.)   
When a local conformal net satisfies these three conditions, we 
say that it is {\sl completely rational}.
We have shown in \cite{KLM} that
if a local conformal net is completely rational,  then we
have only finitely many mutually inequivalent irreducible
representations of the net.  (Such finiteness is often
called rationality.  Our conditions are called complete
rationality since they imply this
rationality.)  Furthermore, a weighted counting
of the number of such representations, where each such 
representation is counted as the square of its statistical
dimension, gives the index in the third condition above \cite{KLM}.
Furthermore, it has been known \cite{FRS} that we have a braiding
on the tensor category of representations of local conformal nets, and 
a notion of non-degeneracy of a braiding has been introduced in
\cite{R1}.  We have also shown in \cite{KLM} that the representation
category of a completely rational local conformal net has a
non-degenerate braiding.  In other words, the tensor category is
modular in the sense of \cite{T}, and we have a unitary
representation of the modular group $SL(2,\Z)$ whose dimension
is the number of unitary equivalence classes of irreducible
representations of the net.

It is difficult in general to show that a certain local conformal net
is completely rational.  F. Xu \cite{X3} proved the finiteness
of the index, the third condition above, for the $SU(N)_k$-nets
constructed by A. Wassermann \cite{W} and thus proved their
complete rationality.  Complete rationality for orbifold and
coset constructions has been also studied by F. Xu
\cite{X4,X5,X6}, where many interesting properties and
examples are presented.  (Also see \cite{L3} in the case of the
coset construction.)  This complete rationality  plays a role
of amenability condition for local conformal nets.

In the classical representation theory, the notion of induced
representation is very important and useful.  Suppose we have
a subgroup $H$ of $G$ and a unitary representation of $H$.
Under a very mild condition, we can produce the {\sl induced
representation} of $G$ from a given representation of $H$.
We would like to perform a similar construction for
representation of local conformal nets.  Suppose we
have a net of subfactors $\A\subset \B$ on the circle.
That is, this is a family of subfactors $\A(I)\subset \B(I)$
parametrized intervals $I\subset S^1$ with a certain natural
set of axioms.  (See \cite{LR,L3}
for a precise set of axioms.)  We assume finiteness of the
index $\A(I)\subset \B(I)$, which is independent of the
interval $I$.  Then, roughly speaking, we can produce a
representation of the larger net $\B$ from a
representation of the smaller net $\A$.  This construction
was introduced by Longo-Rehren \cite{LR} based on an old
suggestion of J. Roberts, and many interesting properties
and examples were studied by F. Xu \cite{X1,X2}.  Further
studies and applications have been made in \cite{BE,BE2,BEK1,BEK2,BEK3}
under the name $\a$-induction.  (This $\a$-induction also generalizes
Ocneanu's graphical construction \cite{O2} for the 
Dynkin diagrams.)
However, we do have three substantial differences between
the induction for group representations and $\a$-induction 
for nets as follows.  First, for groups $H\subset G$, the larger group
has a larger symmetry, but the situation is opposite for nets.
That is, the representation category
of the larger group
$G$ is larger, in an appropriate sense, than that of $H$,
but the representation category of the larger net $\B$ is
{\sl smaller} than that of $\A$.  The second difference is
that $\a$-induction depends on the braiding
structure of the representation category of the smaller
net $\A$.  A braiding always comes in a pair, a positive
braiding and a negative braiding.  (Graphically, they are
represented by over-crossing and under-crossing, respectively.)
We distinguish the two $\a$-inductions by the notations
$\a^+$ and $\a^-$.  The third difference is that
$\a$-induction does not produce a genuine representation
of the larger net and it gives only a ``fake representation'', 
which is often called a {\sl soliton sector}, although it behaves
like a genuine representation in many senses.  
For example, it does have notions
of a tensor product and a statistical dimension, 
but it does not have a braiding.  (See \cite{L1} for a more
general setting of such a tensor category.)
Roughly speaking, the representation
category of the larger net $\B$ is ``too small'' to accept the
image of $\a$-induction.  Actually, the intersection of the set of
irreducible representations arising from a positive $\a$-induction
and that from a negative $\a$-induction precisely gives the
set of irreducible representations of the larger net $\B$.  These
differences produce a new interesting structure that does not
appear in the classical setting of group representations,
as explained in the next section.

We remark that we also have a machinery of restriction of a
representation of a larger net to that of a smaller net, just as in
the case of group representation, but $\a$-induction is
much more useful than the restriction in many cases.
One reason is that
we often have a situation where we have enough information
about the smaller nets and would like to obtain information
on  the larger nets.

\section{Modular invariants and classification of CFT}

Consider a net of subfactors $\A\subset\B$ with finite index as
above, where $\A$ is completely rational.
Choose irreducible representations $\la,\mu$ of the
smaller net $\A$ and apply $\a^\pm$-induction to them.
Let $Z_{\la\mu}$ be $\dim\Hom(\a^+_\la,\a^-_\mu)$.  Roughly
speaking, this is a counting of the intersection of $\a^+_\la$
and $\a^-_\mu$, and such an intersection occurs only in the
genuine representation category of the larger net $\B$, as
explained above.  This number $Z_{\la\mu}$ is independent of
choices of irreducible representations within unitary equivalence
classes, of course, so in this way, we obtain a matrix 
$Z=(Z_{\la\mu})$ of finite size whose rows and columns are
indexed by unitary equivalences classes of irreducible representations
of the net $\A$.  The unitary representation of $SL(2,\Z)$ mentioned
above also has the same set of indices and we have proved 
in \cite[Corollary 5.8]{BEK1} that this matrix $Z$ is in the
commutant of the image of this unitary representation.  The
entries of the matrix $Z$ are obviously non-negative
integers and we have $Z_{00}=1$, where $0$ denotes the vacuum
representation of the net $\A$.  Such a matrix of non-negative
integers in the commutant of a unitary representation of $SL(2,\Z)$,
subject to the normalization $Z_{00}=1$, is called a
{\sl modular invariant} of the unitary representation (or of the
representation category of the net $\A$).  In general, the
unitary representation is not irreducible, so we can have
something non-trivial in the commutant, but  this
unitary representation is often almost irreducible in some sense,
and the commutant is quite small.  So the modular invariant
condition is very restrictive.  In general, the number of
modular invariants for a given unitary representation is always
finite, and this finite number is often very small, such as
1, 2, or 3, in a typical situation.  The most famous classification
list of such modular invariants has been obtained by
Cappelli-Itzykson-Zuber \cite{CIZ} for the representation category
of the $SU(2)_k$-net.  In this case, we have at most three modular
invariants for a given $k$.  In principle, the classification
problem of a given unitary representation of $SL(2,\Z)$ is a
problem in linear algebra and should be ``solvable'', but
the classification can be very hard combinatorially, in practice.
Gannon's series of papers deal with this problem.  (See \cite{G},
for example.)

Our idea \cite{KL} for classification of local conformal nets is using 
such a classification of modular invariants through $\a$-induction.
This is in the spirit of \cite{BMT}.
That is, start with a local conformal net $\B$ with diffeomorphism
covariance.  Then the algebras generated by the projective unitary
representation of $\Diff(S^1)$ produces a subnet $\A\subset \B$.
This is the
Virasoro net having a central charge $c$.  We now assume that
$c<1$ for our classification.  This Virasoro net with $c<1$ falls
in the class of coset nets studied by F. Xu \cite{X4,X5}
from the diagonal inclusion
$SU(2)_{m-1}\subset SU(2)_{m-2}\otimes SU(2)_1$ with
$c=1-6/m(m+1)$, and we know that the Virasoro net $\A$
is completely rational by \cite{L3}.  We can write down the
irreducible representations of the Virasoro nets explicitly
and the number of their unitary equivalence classes is
$m(m-1)/2$.  We can also show that the
inclusion $\A(I)\subset \B(I)$ is irreducible, using strong
additivity of the Virasoro net that holds by \cite{L3},
and then a result in
\cite{ILP} implies that the index of $\A(I)\subset \B(I)$ is
automatically finite.  We can now apply the general machinery
of $\a$-induction and obtain a modular invariant $Z$ from
the inclusion $\A(I)\subset \B(I)$, where the rows and columns
of the matrix $Z$ are indexed by the unitary equivalence
classes of the irreducible representations of the Virasoro
net $\A(I)$.  Note that now the dimension of the
unitary representation of $SL(2,\Z)$ is $m(m-1)/2$.
In this way, we obtain a matrix $Z$ with entries of
non-negative integers and dimension $m(m-1)/2$,
starting from a local conformal net $\B$ with diffeomorphism
covariance and central charge $c=1-6/m(m+1)$, $m=2,3,4,\dots$.
This matrix $Z$ is obviously an invariant of the net $\B$
and subject to a strong constraint of modular invariance.
The modular invariants of this unitary representation of
$SL(2,\Z)$ arising from the Virasoro nets
have already been classified by Cappelli-Itzykson-Zuber
\cite{CIZ} and an explicit classification table is there.
The modular invariants are in a bijective correspondence to
pairs of Dynkin diagrams whose Coxeter numbers have a difference 1.
They are classified into two types.  The type I modular invariants
are labeled with pairs of Dynkin diagrams $A_n$, $D_{2n}$,
$E_6$, $E_8$, while type II modular invariants with
$D_{2n+1}$, $E_7$ (in pairs with the corresponding $A$-diagrams).
(Recall that paragroups of index less than 4 are labeled with
Dynkin diagrams $A_n$, $D_{2n}$, $E_6$, $E_8$ as in \cite{O1,EK}
and appearance of the same diagrams in the context of $SU(2)$
modular invariants is related to locality of the extended nets
as in \cite{BE,BE2,BEK1,BEK2}.)
Locality of the net $\B$ shows that the modular invariant we
obtain from $\B$ as above must be of type I.  So our classification
will be complete if we can show that this map from a net $\B$ to
a modular invariant matrix $Z$ of type I in the classification
list in \cite{CIZ} is bijective.

This is a problem of classifying irreducible
extensions $\B\supset \A$ when the smaller net $\A$ is given.
When a modular invariant $Z$ is given and if it arises from
an extension $\B$, it is easy to see
the decomposition $\bar\iota \iota=\bigoplus_{\la} Z_{0\la}\la$,
where $\iota$ is the inclusion map $\A(I)\subset \B(I)$
and $\la$ varies in a set of representatives of unitary
equivalence classes of irreducible representations of the net $\A$.
If we start with a subfactor $\A(I)\subset \B(I)$, then 
the dual canonical endomorphism $\rho=\bar\iota \iota$ and
intertwiners $V\in\Hom(\id,\rho)$ and $W\in\Hom(\rho,\rho^2)$
satisfy certain relations.  Longo \cite{L2} has shown that
a triple $(\rho,V,W)$ arising from a subfactor can be
axiomatized in terms of simple algebraic relations, and he
has called a triple satisfying these relations a {\sl $Q$-system}.
Thus the classification problem is reduced the classification of
$Q$-systems $(\bigoplus_{\la} Z_{0\la}\la, V, W)$ up to
unitary equivalence in a natural sense.  A $Q$-system arising
from an  extension $\B\supset\A$ should also satisfy the
chiral locality in \cite[Theorem 4.9]{LR}
in the sense of \cite[page 454]{BEK1}.  Such a classification
problem has been studied for  extensions of
the $SU(2)_k$-nets in several different contexts in
\cite{BEK2,KO,O2}.  Using such classification results for the
$SU(2)_k$-nets,
one can also obtain the desired bijectivity for the case
of the Virasoro nets, as in \cite{KL}.  The most subtle
step in the $SU(2)_k$-classification of the $Q$-systems
is uniqueness of the $Q$-system for the $E_8$ modular invariant,
which was first obtained in \cite{KO}, and we use this uniqueness
result for our classification.
Izumi has also made direct computations of this $Q$-system and
proved this uniqueness.  We also use the B\"ockenahauer-Evans
criterion \cite{BE2} for locality of extensions of the Virasoro
nets we construct combinatorially.

The meaning of the labeling in terms of the {\sl pairs} of
Dynkin diagrams in this setting is as follows.  Take an
example of $(A_{10}, E_6)$.  The pair of type $A$ Dynkin
diagrams having the same Coxeter numbers is $(A_{10}, A_{11})$.
This pair of type $A$ diagrams
is realized from the Virasoro net with $c=21/22$.
The representation category of this Virasoro net is 
a direct product of two systems $A_{10}$ and $A_{11}$ divided
by a symmetry of order 2, in an appropriate sense.  (The reason
we have a pair of type $A$ diagrams here was clarified in
the setting of coset constructions by Xu \cite{X4,X5}.  The
number of unitary equivalence classes of irreducible
representations is $10\times 11/2=55$.)
The representation category of the local conformal
net corresponding to the pair $(A_{10}, E_6)$ is a
direct product of two systems $A_{10}$ and $A_3$ divided
by a symmetry of order 2 in the same sense, and the number
of unitary equivalence classes of irreducible
representations is $10\times 3 /2=15$.  In this way, the
meaning of the label $E_6$ is not visible.
However, if we use only the $\a^+$-induction, without $\a^-$-induction,
then we obtain a tensor category of  ``fake'' representations and
the system of the irreducible objects is a direct product of
two systems $A_{10}$ and $E_6$ divided
by a symmetry of order 2 in the same sense above, and the
number of the irreducible objects is $10\times 6 /2=30$.
The $E_6$ system having six irreducible objects does not arise
from a representation category of a net, but it arises from
a tensor category of ``fake'' representations produced with 
$\a^+$-induction.

The final classification result in \cite{KL} is as follows.

\begin{theorem}
The local conformal nets on the circle
with central charge less than $1$ are listed as follows.
\begin{enumerate}
\item Virasoro nets with central charge $c=1-6/m(m+1)$.
\item Their simple current extensions of index $2$.
\item Four exceptional nets labeled with
$(A_{10}, E_6)$, $(E_6, A_{12})$,
$(A_{28}, E_8)$, $(E_8, A_{30})$.
\end{enumerate}
\end{theorem}

Two of the exceptional ones in the above list are realized as
the coset constructions for
$SU(2)_{11}\subset SO(5)_1\otimes SU(2)_1$ and
$SU(2)_{29}\subset (G_2)_1\otimes SU(2)_1$.  Actually, they were
considered by B\"ockenhauer-Evans
\cite[II, Subsection 5.2]{BE} as possible candidates
realizing the corresponding modular invariants in the
Cappelli-Itzykson-Zuber list, but they were unable to prove
that these coset construction indeed produce the desired
nets.  The other two exceptionals appear to be new.

Note that we did not assume the three conditions for
complete rationality in the above, but they
automatically hold.  This is also analogous to the fact
that the (hyperfinite II$_1$) subfactors of index less
than 4 are automatically amenable.  (These subfactors are 
even of finite depth, actually.)   Carpi \cite{C2} has recently
got a classification for local conformal nets
with  $c=1$ arising as compact extensions of
the Virasoro net.

In the proof of the above theorem, we have used the modular invariants
$Z$ arising from  local conformal nets as an invariant, but in
the statement of the theorem, we see
that the representation category does
gives a complete invariant of such local conformal nets.  We also 
have a complete list of such representation categories.  In this
way, our goal of the classification we mentioned at the
beginning of this paper has been achieved for this class of
local conformal nets.  It is expected that this type of
classification theorem holds in a much wider situation, but
proving such a result is beyond the current technique.

We finally
remark why this approach based on $\a$-induction has worked, in
comparison to classification of subfactors.
Here, a local conformal net is regarded as an analogue of a
subfactor.  (They both produce a tensor category through a
``representation theory''.  This gives one analogy.)
We have obtained a Virasoro subnet $\A$ from a given net $\B$.
In an analogy to subfactor theory in this spirit,
the Virasoro subnet corresponds to a
subfactor generated by the downward sequence of the Jones
projections.  That is, for a type II$_1$ subfactor $N\subset M$,
we obtain a commuting square
$$\begin{array}{ccc}
N & \subset & M \\
\cup && \cup \\
P & \subset & Q,
\end{array}$$
where $P$ and $Q$ are generated by the downward sequence of
the Jones projections.  (Recall that a commuting square
is regarded as an embedding of one subfactor into another.)
The Virasoro nets are ``minimal'' in the sense
that it is contained in any diffeomorphism covariant net
and also in the sense of Carpi's work \cite{C}.  Similarly,
a subfactor generated by the Jones projections is ``minimal''
in the sense that it is contained in any subfactor.
But the difference is that the Virasoro net is canonically
defined while the above subfactor $P\subset Q$ is not
canonically determined in $N\subset M$.  In this way, we
have  a control over the index $\A \subset \B$ while we
have no control over the index $P\subset N$.  This explains
why our approach has worked while a similar approach to 
classification of subfactors with index less than 4 seems
much more difficult (although we have a by far more general
classification result for subfactors in \cite{P}).

\medskip
\noindent{\bf Acknowledgments.}
We are grateful to S. Carpi for pointing out inaccuracies in
the draft of this paper and  R. Longo for detailed comments on it.
We also thank J. Fuchs, K.-H. Rehren, C. Schweigert, and the anonymous
referee for useful comments.

\bigskip
\bigskip
{\bf Added in Proof:}
After the submission of this paper, the following related two
preprints appeared.

S. K\"oster, 
{\it Local nature of coset models}, preprint 2003,
math-ph/0303054.

F. Xu, 
{\it Strong additivity and conformal nets}, preprint 2003,
math.QA/0303266.

In the former, K\"oster identified one of our two ``new'' nets,
$(A_{10}, E_6)$, with the two cosets $SU(9)_2\subset (E_8)_2$ and
$(E_8)_3\subset (E_8)_2\otimes (E_8)_1$, assuming that the local
conformal nets $(E_8)_k$ have the expected WZW-fusion
rules.  The remaining one, $(A_{28}, E_8)$, does not seem to be
a coset, and it appears to be a genuine new example.

In the latter, Xu classified the local conformal nets with $c=1$
with an extra assumption.
\end{document}